\documentclass[12pt,a4paper]{amsart}
\usepackage{amsmath,amscd}
\usepackage{latexsym,amsmath,amssymb,mathrsfs,setspace,enumerate,tikz}
\usepackage[T1]{fontenc}
\usepackage[utf8]{inputenc}
\usepackage{lmodern}
\usepackage{amssymb}
\usepackage{float}
\usetikzlibrary{arrows, automata}
\usepackage{longtable, pdflscape}
\usepackage{pgf, tikz}
\usetikzlibrary{arrows, automata}
\newtheorem{lemma}{Lemma}
\newtheorem{prop}{Proposition}
\newtheorem{thm}{Theorem}
\newtheorem{rem}{Remark}
\newtheorem{exa}{Example}
\newtheorem{cor}{Corollary}
\newtheorem{dfn}{Definition}

\newcommand{\bl}{\begin{lemma}}
	\newcommand{\el}{\end{lemma}}
\newcommand{\bt}{\begin{thm}}
	\newcommand{\et}{\end{thm}}
\newcommand{\bc}{\begin{cor}}
	\newcommand{\ec}{\end{cor}}
\newcommand{\bd}{\begin{dfn}}
	\newcommand{\ed}{\end{dfn}}
\newcommand{\bp}{\begin{proof}}
	\newcommand{\ep}{\end{proof}}
\newcommand{\br}{\begin{rem}}
	\newcommand{\er}{\end{rem}}
\newcommand{\be}{\begin{exa}}
	\newcommand{\ee}{\end{exa}}
\newcommand{\bca}{\begin{case}}
	\newcommand{\eca}{\end{case}}
\newcommand{\bcl}{\begin{cla}}
	\newcommand{\ecl}{\end{cla}}
\newcommand{\bil}{\begin{ill}}
	\newcommand{\eil}{\end{ill}}
\newcommand{\bpr}{\begin{prop}}
	\newcommand{\epr}{\end{prop}}
\def \i{\indent}

\def \mc{\mathcal}
\def \tt{\textit}
\date{08.07.19}
\title{\textbf {Ideal categories of rings and \\the ring of cones}}
\author {Sreejamol P.R.$^1$ and P.G. Romeo$^2$}
\address{$^{1}$  Assistant Professor, Department of Mathematics, SNM College, Maliankara, Ernakulam, Kerala, India.,$^2$ Professor, Dept. of Mathematics, Cochin University of Science and Technology, Kochi, Kerala, India.}
\email{$sreejasooraj@gmail.com,\, romeo_-parackal@yahoo.com $}

\begin{document}
\begin{abstract}
In this paper we describe the ideal category of a ring $R$ as preadditive proper category. Further it is also shown that the cones in this category is a ring with appropriate addition and multiplication. 
\end{abstract}
\maketitle
\section{Introduction}
  This paper is related to the study of ideals of a ring. Recall that a regular semigroup $S$ is a semigroup in which for every $a\in S$, there exists $b\in S$ such that $aba = a$. In 1994, K.S.S.Nambooripad introduced a special type of categories called normal categories to describe the ideal structure of a regular semigroup. The categories of principal left (right) ideals of a regular semigroup are normal. Later on, this approach was extended to arbitrary semigroups by defining the ideal categories of a semigroup as proper categories and it is shown that the set of all proper cones in such a category is a semigroup. Also the dual of the proper category is obtained(cf.\cite{PGS}).\\
  \i In the 1930's von Neumann introduced regular rings in his work on continuous geometry. Regular ring (von Neumann regular) is a ring in which the multiplicative part is a regular semigroup. In 
  \cite{Sunny}  the principal left[right] ideals of a regular ring are described as preadditive normal categories. In this paper we extend the categorical approach to the study of the structure theory of rings. We show that the principal left[right] ideals of a ring are preadditive proper categories, further we introduce $RR-$proper categories and it is shown that the set of all proper cones in such a category becomes a ring. Note that all rings in this paper are associative with unit.

\section{Prelimanires}
In the following we briefly recall some definitions and basic results regarding categories needed in the sequel and for more details reader is referred to S. Maclane (cf.\cite{Mac}) and K.S.S.Nambooripad (cf.\cite{KSS}). The categories we are considering are all small categories - a small category is a category in which the class of objects and class of morphisms are both sets. Let $\mc C$ be any category, then $\tt{v}\mc C$ denotes the set of objects of  $\mc C$. For $a, b \in \tt v\mc C$, the set consists of all morphisms of the category with domain $a$ and codomain $b$ is called the 
\tt{hom-set} and is denoted as $\mc C(a;b)$ or $hom(a;b)$. If $f: a\rightarrow b$ and $g: b\rightarrow c$ then the composition is $f\cdot g: a\rightarrow c$. A morphism $f$ in a category $\mc C$ is a monomorphism if for $g,h\in \mc C$, $gf = hf$ implies $g = h$; that is $f$ is a monomorphism if it is right cancelable. Dually a morphism $f\in \mc C$ is an epimorphism if $f$ is left cancelable. \\
\i  An object $a$ is terminal in $\mc C$ if for each object $b$ there is exactly one arrow $b\rightarrow a$. An object $c$ is initial object if to each object $b$ there is exactly one arrow $c\rightarrow b.$ A \tt{zero object} or \tt{null object} $z$ in  $\mc C$ is an object which is both initial and terminal. For any two objects $a$ and $b$ the unique arrows $a\rightarrow z$ and $z\rightarrow b$ have a composite $O_a^b: a\rightarrow b$ called the \tt{zero arrow} from $a\rightarrow b$.

\bd(cf.\cite{Mac})
	A preadditive category(or \tt {Ab}-category) is a category satisfying:
\begin{enumerate}
\item each hom-set is an additive abelian group,
 \item composition of arrows is bilinear relative to this addition,
 \item category has zero object. 
\end{enumerate}
\ed
\i A preorder $\mc P$ is a category such that, for any $p,p'\in \tt v\mc P$; the hom-set $\mc P(p,p')$ contains at most one morphism. In this case, the relation $\subseteq$ on the class $\tt v\mc P$ defined by $p\subseteq p' \Leftrightarrow \mc P(p,p') \neq \phi $ is a quasi-order on $\tt v\mc P$ . In a preorder, p and p' are isomorphic if and only if $\mc P(p,p') \neq \phi \neq \mc P(p',p)$.  Therefore $p\subseteq p'$ is a partial order if and only if $\mc P$ does not contain any nontrivial isomorphisms. Equivalently, the only isomorphisms of $\mc P$ are identity morphisms, and in this case $\mc P$ is said to be a $\tt{strict preorder}.$ \\
\i Let $\mc C$ be a category and $\mc P$ be a subcategory of $\mc C$. Then $(\mc C,\mc P)$ is called a \tt{category with subobjects} if the following hold:
\begin{itemize}
  \item [(1)]$\mc P$ is a strict preorder with $\tt v\mc P = \tt v\mc C$.
  \item [(2)] Every $f \in \mc  P$ is a monomorphism in $\mc C$.
  \item [(3)] If $f, g \in \mc P$ and if $f = hg$ for some $h \in \mc C$, then $h \in \mc P$.
\end{itemize}
In a category with subobjects, if $f: a\rightarrow b $ is a morphism in $\mc P$, then $f$ is an \tt{inclusion} and we denote this inclusion by 
$j(a,b)$. If there is a morphism $e: b\rightarrow a$ such that 
$j(a,b)e = I_a$, then $e$ is called a retraction from $b\rightarrow a$ and is denoted by $e(b,a).$ In case a retraction from b to a exists, then we say that the inclusion $j(a,b): a\rightarrow b$ splits.\\
 \i A morphism $f\in \mc C$, where $\mc C$ be a category with subobjects  has factorization if $f = p\cdot m$ where $p$ is an epimorphism and $m$ is an embedding. A category $\mc C$ is said to have the \tt {factorization property} if every morphism of $\mc C$ has a factorization. Thus, if $\mc C$ has the factorization property, then any morphism $f$ in $\mc C$ has atleast one factorization of the form $f = qj$, where $q$ is an epimorphism and $j$ is an inclusion and such factorizations are called \tt{canonical factorizations}. A \tt{normal factorization} of a morphism $f$ in $\mc C$ is a factorization of the form $f = euj$ where $e$ is a retraction, $u$ is an isomorphism and $j$ is an inclusion.\\
 \i A morphism $f$ in a category with subobjects is said to have an \tt{image} if it has a \tt{canonical(epi-mono)factorization} $f = xj$, where $x$ is an epimorphism and j is an inclusion with the property that whenever $f = yj'$ is any other canonical factorization, then there exists an inclusion $j''$ such that $y = xj''$(\cite{KSS}). A category is said to have \tt{images} if every morphism in $\mc C$ has an image. In this case, the codomain of $x$ is said to be the \tt{image} of $f$. When the morphism $f$ has an image we denote the unique canonical factorization of $f$ by $f = f^oj_f$, where $f^o$ is the \tt{unique epimorphic component} and $j_f$ is the inclusion of $f$.\\
\indent Let $\mc C$ be a category with subobjects, images, every morphism in $\mc C$ has normal factorizations in which the inclusion splits and $d \in \tt v\mc C$. A cone from $\tt v\mc C$ to the vertex $d$ is a map $\gamma: \tt v\mc C \rightarrow d$ such that
\begin{enumerate}
\item $\gamma(c) \in \mc C(c,d)$ for all $c \in \tt v\mc C.$
\item  If $c'\subseteq c$ then $j(c',c) \gamma(c) = \gamma(c')$.
\end{enumerate}
The cone $\gamma$ is called a \tt{normal cone} if there exists an $a\in \tt v\mc C$ such that $\gamma(a)$ is an isomorphism. The vertex $d$ of the cone $\gamma$ is usually denoted as 
$c_\gamma$. The set of all normal cones in $\mc C$ is denoted by $\mc{TC}$ which is a regular semigroup. Corresponding to the normal cone $\gamma$, an $M-$set is defined as 
$$M_\gamma = \{c\in \tt{v}\mc C: \gamma(c) \text{ is isomorphism}\}.$$ 
\bd
A normal category is a pair $(\mc C,\mc P)$ satisfying the following:
\begin{enumerate}
\item $(\mc C,\mc P)$ is a category with subobjects
\item every inclusion in $\mc C$ splits
\item Any morphism in $\mc C$ has a normal factorization
\item for each $a\in \tt v\mc C$ there is a normal cone $\gamma$ with vertex $a$ and $\gamma(a) = I_a.$
\end{enumerate}
\ed
\bl[\cite{KSS}]
 If $\gamma$ is a normal cone with vertex $c_\gamma$ and $f: c_\gamma\rightarrow a$ is an epimorphism for $a\in \tt v\mc C$, then the map 
 $\gamma \star f$ defined by 
$$(\gamma\star f)(a) = \gamma(a)\cdot f$$
from $\tt v\mc C \rightarrow \mc C$is a normal cone with vertex $d$.
\el

\bt
Let $\mc C$ be a normal category. Then set of all normal cones $\mathcal{T}\mc C$ with respect to the binary operation defined by 
$$\gamma\cdot\eta = \gamma\star\eta(c_\gamma)^o$$ 
for all $\gamma,\eta \in \mathcal{T}\mc {C}$ is a regular semigroup. 
\et

\section{Proper categories}
In the following we proceed to describe certain generalization of normal categories in which the cones at each vertex need not always be normal but only satisfies some weaker conditions and we call proper categories. 
\bd
Let $\mc C$ be a category with subobjects, every inclusion splits and every morphism has canonical factorization. A proper cone $\gamma$ in $\mc C$ is a cone with vertex $d$ such that there exists at least one $c\in \tt{v}\mc C$ with $\gamma(c): c\rightarrow d$ is an epimorphism (i.e., $\gamma(c) = \gamma(c)^o$). 
\ed 

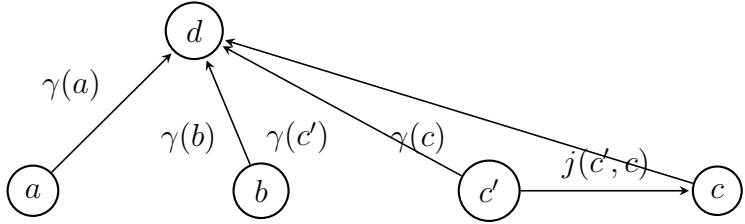
\begin{tikzpicture}[
> = stealth, % arrow head style
shorten > =1pt, % don't touch arrow head to node
auto,
node distance = 3cm, % distance between nodes
semithick % line style
]

\tikzstyle{every state}=[
draw = black,
thick,
fill = white,
minimum size = 4mm
]

\node[state] (a) {$a$};
\node[state] (d) [above right of=a] {$d$};
\node[state] (b) [right of=a] {$b$};
\node[state] (c') [right of=b] {$c'$};
\node[state] (c) [right of=c'] {$c$};
\path[->] (a) edge node {$\gamma(a)$} (d);

\path[->] (b) edge node {$\gamma(b)$} (d);
\path[->] (c) edge node {$\gamma(c)$} (d);
\path[->] (c') edge node {$j(c',c)$} (c);
\path[->] (c') edge node {$\gamma(c')$} (d);

\end{tikzpicture}\\
\i  The set of all proper cones in category $\mc C$ is denoted by $\mc {PC}$ and for  $\gamma \in \mc {PC}$, we denote by $c_\gamma$ the vertex of $\gamma$ and by $M_\gamma$, the $M-$set defined by
$$M_\gamma = \{c\in \tt{v}\mc C: \gamma(c) \text{ is epimorphism}\}$$
A cone $\gamma$ in $\mc C$ is proper or normal cone according as $M_\gamma \neq \phi$. 
Clearly every normal cone is proper as well.
\bd
A small category $\mc C$ with subobjects is called \tt{proper category} if it satisfies the following:
\begin{enumerate}
	\item every inclusion in $\mc C$ splits,
	\item any morphism $f\in \mc C$ has canonical factorization and
	\item each object of $\mc C$ is a vertex of a proper cone $\gamma\in \mc {PC}$. 
\end{enumerate} 
\ed
\bd
A preadditive category which is also a proper category is termed as a \tt{preadditive proper category}.\\
\i A preadditive proper category $\mc C$ satisfying the following conditions: 
\begin{enumerate}
\item the object set $\tt  v\mc C$ with partial order induced by subobject relation is a relatively complemented lattice and 
\item every subset of the object set which has an upper bound contains a unique maximal element 
\end{enumerate} 
is called an \tt{$RR-$proper category} and these conditions are termed as $RR-$conditions.
\ed
\i Now we proceed to show that the set of proper cones $\mc{PC}$ in an $RR-$ proper category $\mc C$, is a ring.

\bpr
Let $\mc C$ be an $RR-$proper category and $\gamma$ be a proper cone in $\mc C$. For $f\in \mc C(c_\gamma,c)$, $\gamma\star f^o$ is a proper cone with vertex $im f$ and the components are given by $(\gamma\star f^o)(a) = \gamma(a)f^o $ for all $a\in \tt v\mc C$ such that for all composable pair of morphisms $f, g\in \mc C$ with $dom f = c_\gamma$
$$\gamma\star (fg)^o = (\gamma\star f^o)\star (j_{c_1}^cg)^o$$
where $c_1 = im f = imf^o$.
\epr 
\bp
 For $c\in M_\gamma$, $\gamma(c)$ is an epimorphism and so is $\gamma(c)\cdot f^o$. Hence $c\in M_{\gamma\star f^o}$ and hence $\gamma\star f^o$ is a proper cone. \\
Now $f$ and $g$ are as in the statement. Let $\eta = \gamma\star f^o$ so that $c_\eta = im f = c_1$. Then we have
\begin{align*}
(\gamma\star (fg)^o)(a) &= \gamma(a)f^o(j_{c_1}^cg)^o\\
& = \eta\star(j_{c_1}^cg)^o(a)\\
& =  (\gamma\star f^o)\star (j_{c_1}^cg)^o(a)
\end{align*}
for all $a\in \tt{v}\mc C.$
\ep
\bt
Let $\mc C$ be a proper category. Then $\mc {PC}$ the set of all proper cones in $\mc C$, is a semigroup with respect to the binary operation defined by $$\gamma\cdot\eta = \gamma\star\eta(c_\gamma)^o$$ 
for all $\gamma,\eta \in \mc {PC}$.
\et 
\bp
For $\gamma, \eta \in \mc {PC}$,  $\gamma\cdot\eta$ is a proper cone with vertex $c_{\gamma\cdot\eta} = im\text{}\eta(c_\gamma)$. To show that the binary operation defined is associative, let $\alpha,\beta,\gamma \in \mc {PC}$ and for $c\in \tt{v}\mc C,$ 
\begin{align*}
(\alpha(\beta\gamma))(c) &= \alpha(c)((\beta\gamma)(c_\alpha))^o\\
& = \alpha(c)\huge(\beta(c_\alpha)((\gamma(c_\beta))^o\huge)^o\\
& = \alpha(c)[(\beta(c_\alpha))^o j_{Im \beta(c_\alpha)}^{c_\beta}(\gamma(c_\beta))^o]^o\\
& = \alpha(c)((\beta(c_\alpha))^o(\gamma(c_{\alpha\beta}))^o\\
((\alpha\beta)\gamma)(c) &= (\alpha\beta)(c) (\gamma(c_{\alpha\beta}))^o\\
&= \alpha(c)((\beta(c_\alpha))^o(\gamma(c_{\alpha\beta}))^o
\end{align*}
Thus $\alpha(\beta\gamma) = (\alpha\beta)\gamma$ and hence $\mc {PC}$ is a semigroup.
\ep
\bpr
$\gamma\in \mc {PC}$ is an idempotent proper cone if and only if $\gamma(c_\gamma) = I_{c_\gamma}$.
\epr
\bp
Suppose $\gamma$ is an idempotent proper cone and let $c\in M_\gamma$. Then $(\gamma\cdot\gamma)(c) = \gamma(c)$ implies $\gamma(c)(\gamma(c_\gamma))^o = \gamma(c)$. Since $\gamma(c)$ is an epimorphism $(\gamma(c_\gamma))^o = I_{c_\gamma}$. Clearly $\gamma(c_\gamma)\in \mc C(c_\gamma,c_\gamma)$ and so $\gamma(c_\gamma) = I_{c_\gamma}$. Conversely, if $\gamma(c_\gamma) = I_{c_\gamma}$, then for every $a\in \tt{v}\mc C$, $(\gamma\cdot\gamma)(a) = \gamma(a)(\gamma(c_\gamma))^o = \gamma(a)I_{c_\gamma}= \gamma(a)$. Hence $\gamma$ is an idempotent proper cone.\\
\i $E(\mc {PC})$ denotes the set of all idempotent proper cones in $\mc C$.
\ep

\bpr[\cite{Sunny}]
 Let $\gamma$ be a cone with vertex $d$ in an $RR-$proper category $\mc C$. Let $X = \{im\gamma(a): a\in \tt v\mc C$\} then there exists a unique maximum element $d_0\leq d$ in $X$.  
\bp 
Since $RR-$proper category $\mc C$ satisfies $RR-$conditions, $X$ has a maximal element say $d_0$. To prove that it is unique let $d_0$ and $d_1$ be two maximal elements of $X$. Then there exists $b, c\in \tt v\mc C$ such that $d_0 = im \gamma(b)$ and $d_1 = im \gamma(c)$. By $RR-$conditions $b\vee c\in \tt v\mc C$. Let $d = b\vee c.$ Then $b\subseteq d$ and $c\subseteq d.$ Then 
$$j(b,d)\gamma(d) = \gamma(b)\text{ and }j(c,d)\gamma(d) = \gamma(c).$$
Obviously
$$im j(b,d)\gamma(d) \subseteq im \gamma(d)\text{ and } im j(c,d)\gamma(d) \subseteq im \gamma(d).$$
Then $im \gamma(b) \subseteq im \gamma(d)$ and $im \gamma(c) \subseteq im \gamma(d)$ so that $d_0 \subseteq im \gamma(d)$ and $d_1 \subseteq im \gamma(d)$. Since $d_0$ and $d_1$ are maximal elements in $X$, $d_0 = im \gamma(d)$ and $d_1 = im \gamma(d)$. So $d_0 = d_1$. Also $im \gamma(a)\leq d_0$ for all $a\in \tt v\mc C$.  
\ep 
\epr   
\bl 
If $\gamma$ is a cone in an $RR-$proper category $\mc C$ with vertex $d$ and $d_0 = max\{im\gamma(a)|a\in \tt v\mc C\}$ then for every retraction $e: d\rightarrow d_0$, the cone $\gamma^*$ with vertex $d_0$, defined by;
$$\gamma^*(a) = \gamma(a)e(d,d_0), \text{ for all }a\in \tt v\mc C$$
is a proper cone.
\el
\bp 
For all $a\in \tt v\mc C, \gamma^*(a): a\rightarrow d_0$. Since $e: d\rightarrow d_0$ is a retraction, $d_0\subseteq d$. Let $d_0 = im \gamma(d_1)$ where $d_1\in \tt v\mc C$.
\begin{align*}
\text{ For }a\subseteq b, j_a^b \gamma^*(b) &= j_a^b \gamma(b)e(d,d_0) \\
&= \gamma(a)e(d,d_0)\text{ } [\text{since } \gamma \text{ is a cone }]\\
&= \gamma^*(a)
\end{align*}
hence $\gamma^*$ is a cone. Now to prove that $\gamma^*$ is a proper cone, it is sufficient to prove that at least one component of $\gamma^*$ is an epimorphism.\\
We have $im \gamma(d_1) = d_0$ and $\gamma(d_1) = qj$ is the canonical factorization, where $q$ is the epimorphism and $j: d_0\subseteq d$, inclusion, then 
\begin{align*}
\gamma^*(d_1) &= \gamma(d_1)e(d,d_0)\\
&= qj(d_0,d)e(d,d_0)\\
&= q
\end{align*}
which is an epimorphism and hence $\gamma^*$ is a proper cone.
\ep 
\br  If $\gamma$ is a proper cone, then $\gamma^* = \gamma.$
\er
\bl
Let $\gamma^*$ be a proper cone in the $RR-$proper category $\mc C$ as defined in Lemma 2, then the epimorphic component of $\gamma^\ast$ is $\gamma^\ast$ i.e., $((\gamma^\ast)(a))^o = \gamma^\ast(a),\forall a\in \tt v\mc C.$ 
\el
%\bp 
%$(\gamma^\ast)^o(a) = [\gamma(a)e(c_\gamma,d_0)]^o$ where  $d_0 = max\{im \gamma(a)|a\in \tt v\mc C\}$. Since $e(c_\gamma,d_0)$ is the retraction, $[\gamma(a)e(c_\gamma,d_0)]^o = \gamma(a)e(c_\gamma,d_0)$. Hence $(\gamma^\ast)^o(a) = \gamma^\ast(a).$
%\ep  
\bl
If $\gamma$, $\beta$ are two proper cones in an $RR-$proper category $\mc C$ and $\gamma^*$ is defined as in Lemma 2, then
$$(\gamma\cdot \beta)^* = \gamma\cdot \beta^*$$
\el
\bp 
For all $a\in \tt v\mc C$, and let $d_0 = max\{im (\gamma\cdot \beta)(a)|a\in \tt v\mc C\}$
\begin{align*} 
(\gamma\cdot \beta)^*(a) &= (\gamma\cdot \beta)(a)\cdot e(c_\beta,d_0)\\
& = \gamma(a)\cdot [(\beta(c_\gamma))^o\cdot e(c_\beta,d_0)]\\
& = \gamma(a)\cdot [\beta^*(c_\gamma)]^o\\
& = \gamma(a)\cdot [\beta^*(c_\gamma)]\text{ }[\text{ by Lemma 3}]\\
& = (\gamma\cdot \beta^*)(a)
\end{align*} 
Hence the proof.
\ep 
\bl[\cite{Sunny}]
Let $\mc C$ be an $RR-$proper category and $\mc {PC}$ the semigroup of proper cones in $\mc C$. For $\gamma, \delta \in \mc {PC}$, with vertices $c_\gamma = c$ and $c_\delta = d$  and for all $a\in \tt v\mc C$,
$$ (\gamma\oplus\delta)(a) = \gamma(a)j(c,c\vee d) + \delta(a)j(d,c\vee d)$$  
Then $ \gamma\oplus\delta$ is a \tt{cone} with vertex $c\vee d$.
\el
\bp

$\gamma(a)j(c,c\vee d)$ and $\delta(a)j(d,c\vee d)$ are in $hom(a,c\vee d)$.
Let $a\subset b$. Then 
\begin{align*}
j(a,b)(\gamma\oplus\delta)(b) &= j(a,b)[\gamma(b)j(c,c\vee d) + \delta(b)j(d,c\vee d)]\\
&= j(a,b)\gamma(b)j(c,c\vee d) + j(a,b)\delta(b)j(d,c\vee d)]\\
&= \gamma(a)j(c,c\vee d) + \delta(a)j(d,c\vee d)\\
&= (\gamma\oplus\delta)(a)
\end{align*}
Thus $ \gamma\oplus\delta$ is a \tt{cone} with vertex $c\vee d$.
\ep 
\bc
If $c_\gamma = c_\delta$, then the inclusions $j(c_\gamma,c_\gamma\vee c_\delta)$ and $j(c_\delta,c_\gamma\vee c_\delta)$ become identity maps and $ (\gamma\oplus\delta)(a) = \gamma(a) + \delta(a)$.
\ec 
\bd 
Let $\mc C$ be an $RR-$proper category and $\mc {PC}$ the semigroup of proper cones in $\mc C$. For $\gamma, \delta \in \mc {PC}$, with vertices $c_\gamma = c$ and $c_\delta = d$, define $ \gamma+\delta = (\gamma\oplus\delta)^\ast$ where $(\gamma\oplus\delta)^\ast(a) = (\gamma\oplus\delta)(a) e(c\vee d,d_0)$ for all $a\in \tt v\mc C$. Then $ \gamma+\delta$ is a proper cone  with vertex $d_0$, where $d_0 = max\{im (\gamma\oplus\delta)(a), /a\in \tt v\mc C\}$.
\ed
\bl
For an $RR-$proper category $\mc C$, the set of all proper cones $\mc {PC}$ with  the addition defined by $ \gamma+\delta = (\gamma\oplus\delta)^\ast$, for $\gamma, \delta \in \mc {PC}$ is an additive abelian group.
\el
\bp
For $\gamma, \delta, \rho \in \mc {PC}$, $ \gamma+\delta$ is a proper cone in $\mc {PC}$. Since $\mc C$ is a $RR-$proper category, it is preadditive and hence each homset in $\mc C$ is an additive abelian group.\\
Now let $\gamma + \delta = \eta$, $\delta + \rho = \tau, (c_\gamma\vee c_\delta)\vee c_\rho = c_\gamma\vee (c_\delta\vee c_\rho) = d$, then  
\begin{align*}
((\gamma + \delta) + \rho)(a) &= (\eta + \rho)(a)\\
&= (\eta \oplus \rho)^*(a)\\
&= [\eta(a)j(c_\eta,d) + \rho(a)j(c_\beta,d)]^*\\
&= [(\gamma(a)j(c_\gamma,c_\gamma\vee c_\delta) \\
&+\delta(a)j(c_\delta,c_\gamma\vee c_\delta))e(c_\gamma\vee c_\delta,c_\eta)j(c_\eta,d) + \rho(a)j(c_\rho,d)]^*\\
&= [(\gamma(a)j(c_\gamma,d) + \delta(a)j(c_\delta,d)) + \rho(a)j(c_\rho,d)]^*\\
&=  [\gamma(a)j(c_\gamma,d) + (\delta(a)j(c_\delta,d) + \rho(a)j(c_\rho,d))]^*
\end{align*}
and 
\begin{align*}
(\gamma + (\delta + \rho))(a) &= (\gamma + \tau)(a)\\
&= (\gamma \oplus \tau)^*(a)\\
&= [\gamma(a)j(c_\gamma,d) + \tau(a)j(c_\tau,d)]^*\\
&= [\gamma(a)j(c_\gamma,d) + (\delta(a)j(c_\delta,c_\delta\vee c_\rho) \\
&+ \rho(a)j(c_\rho,c_\delta\vee c_\rho))e(c_\delta\vee c_\rho,c_\tau)j(c_\tau,d) ]^*\\
&= [\gamma(a)j(c_\gamma,d) + (\delta(a)j(c_\delta,d) + \rho(a)j(c_\rho,d))]^* 
\end{align*}
hence the addition is associative.\\
Let $0$ be the zero object in $\mc C$ and $\gamma_0$ be the cone with vertex $0$, where $\gamma_0(a) = 0$ for all $a\in \tt v\mc C$, then $\gamma_0(a)$ is the unique morphism from $a$ to $0$ and $\gamma_0$ is a proper cone in $\mc {PC}$. For every $\gamma \in \mc {PC}$ and for all $a \in \tt v\mc C$, let $d= max\{im(\gamma \oplus \gamma_0)(a)|a\in \mc C\} = c_\gamma$, since $\gamma$ is a proper cone.
\begin{align*}
(\gamma + \gamma_0)(a) &= (\gamma \oplus \gamma_0)^*(a)\\
&= [\gamma(a)j(c_\gamma,c_\gamma) + \gamma_0(a)j(c_{\gamma_0},c_\gamma)]e(c_\gamma,d)\\
& = [\gamma(a) + \gamma_0(a)j(0,c_\gamma)]e(c_\gamma,d)\\
& = \gamma(a)
\end{align*}
thus  $\gamma + \gamma_0 = \gamma$. Similarly  $\gamma_0 + \gamma =  \gamma$. Hence $\gamma_0$ is the identity element in $\mc {PC}$.\\
\i For $\gamma \in \mc {PC}$, define $-\gamma$ by $-\gamma(a) = -(\gamma(a))$. 
Clearly $-\gamma$ is a proper cone in $\mc {PC}$ and $c_{-\gamma} = c_\gamma$.
\begin{align*}
(\gamma + -(\gamma))(a) &= [\gamma(a)j(c_\gamma,c_\gamma\vee c_{-\gamma}) + (-\gamma)(a)j(c_{-\gamma},c_\gamma\vee c_{-\gamma})]^*\\
&= [\gamma(a)j(c_\gamma,c_\gamma) + (-\gamma)(a)j(c_\gamma,c_\gamma)]^*\\
&= [\gamma(a)I_{c_\gamma} + (-\gamma)(a)I_{c_\gamma}]^*\\
&= [\gamma(a) + (-\gamma)(a)]^* \\
&= [\gamma_0(a)]^* = \gamma_0(a) 
\end{align*}
i.e., $\gamma + (-\gamma) = \gamma_0$. 
\begin{align*}
(\gamma + \delta)(a) &= [(\gamma \oplus \delta)(a)]^*\\
&= [\gamma(a)j(c_\gamma,c_\gamma\vee c_\delta) + \delta(a)j(c_\delta,c_\gamma\vee c_\delta)]^*\\
&= [\delta(a)j(c_\delta,c_\gamma\vee c_\delta) + \gamma(a)j(c_\gamma,c_\gamma\vee c_\delta)]^*\\
& = [(\delta \oplus \gamma )(a)]^*\\
&= (\delta + \gamma )(a)
\end{align*}	
Thus $\mc {PC}$ is an additive abelian group.
\ep  
\bt
For an $RR-$proper category $\mc C$, the set of all proper cones $\mc{PC}$ is a ring.
\et
\bp 
For $RR-$proper category $\mc C$, the set of all proper cones $\mc {PC}$ is a semigroup with respect to the multiplication, 
$$\gamma\cdot\beta = \gamma\star\beta(c_\gamma)^o$$
where $\gamma,\beta \in \mc{PC}$. Since $\mc{PC}$ is an additive abelian group it is enough to prove that the multiplication  distributes over addition. As $RR-$proper category is preadditive, composition of morphisms is distributive over addition.\\
For $\gamma, \delta, \rho \in \mc{PC}$, and for all $a\in \tt v\mc C$, let $c_\rho = c$ and $c_\gamma\vee c_\delta = d$
\begin{align*}
[\rho\cdot (\gamma + \delta)](a) &= \rho(a)\cdot [(\gamma + \delta)(c)]^o \\
& = \rho(a)\cdot [(\gamma \oplus \delta)^*(c)]^o \\
& = \rho(a)\cdot ([(\gamma(c)j(c_\gamma,d) + \delta(c)j(c_\delta,d))]^*)^o\\
& = \rho(a)\cdot [(\gamma(c)j(c_\gamma,d) + \delta(c)j(c_\delta,d))]^*\text{ }[\text{by Lemma 3 }]
\end{align*}
\begin{align*}
[(\rho\cdot \gamma) + (\rho\cdot \delta)](a) &= [(\rho\cdot \gamma) \oplus (\rho\cdot \delta)]^*(a)\\
& = [\rho(a)\cdot (\gamma(c))^oj(im \gamma(c),d) +  \rho(a)\cdot (\delta(c))^oj(im \delta(c),d)]^*\\
& = [\rho(a)\cdot (\gamma(c))^oj(im \gamma(c),c_\gamma)j(c_\gamma,d) \\
& +  \rho(a)\cdot (\delta(c))^oj(im \delta(c),c_\delta)j(c_\delta,d)]^*\\ 
& = [\rho(a)\cdot (\gamma(c)j(c_\gamma,d) + \delta(c)j(c_\delta,d))]^*\\ 
& = \rho(a)\cdot [(\gamma(c)j(c_\gamma,d) + \delta(c)j(c_\delta,d))]^*\text{ }[\text{by Lemma 4 }]
\end{align*} 
Hence $\rho\cdot (\gamma + \delta) = (\rho\cdot \gamma) + (\rho\cdot \delta)$.
Similarly $(\gamma + \delta)\cdot \rho = \gamma\cdot \rho + \delta\cdot \rho$, thus $\mc{PC}$ is a ring.
\ep

\section{Ideal categories of rings}
In this section we describe the principal left[right] ideal categories of a ring and we show that they form $RR-$proper categories. Let $R$ be a ring. Then principal left [right] ideal of a ring $R$ generated by an element $a$ of $R$ is $(a)_l = Ra$ [$(a)_r = aR$].\\
\i Since the multiplicative part of a ring is semigroup, we can show that the principal left[right] ideals of a ring $R$ as objects and morphisms right[left] translations form \tt{proper category} $\mathbb L(R)[\mathbb R(R)]$.

\bl Let $R$ be a ring. Then $\mathbb L(R)$, the set of principal left ideals of $R$, is a category whose objects and morphisms are as defined below.
$$\tt v\mathbb L(R) = \{Ra: a\in R\} \text{ and for }a,b \in R $$
$$\mathbb L(R)(Ra,Rb) = \{\rho(a,s,b): x\rightarrow xs\text{ with }as\in Rb; \text{ for all }x\in Ra\}.$$ Then $\mathbb L(R)$ is a category.
\el
\bp
For all $a,b \in R$ and $x\in Ra$, $\rho(a,s,b): x\rightarrow xs\text{ with }as\in Rb$ is the map $\rho_s|_{Ra} \text{ such that } s\in R \text{ and } as\in Rb$. The composition in $\mathbb L(R)$ is given by the rule \\
$\rho(a,s,b)\cdot\rho(c,t,d) =
\begin{cases}
\rho(a,st,d),\text{ if } Rb = Rc  \\

undefined, \text{ otherwise}.
\end{cases}$\\
is associative whenever the composition is defined and $\rho(a,1,a) = I_{Ra}: Ra\rightarrow Ra$ for all $a\in R$ is the identity morphism and hence $\mathbb L(R)$ is a category.
\ep 
Dually $\mathbb R(R)$ is also a category with objects principal right ideals and  morphisms left translations,
$$\tt v\mathbb R(R) = \{aR: a\in R\}\text{ and}$$ 
$$\mathbb R(R)(aR,bR) = \{\lambda(a,s,b): x\rightarrow sx\text{ with } x\in aR\text{ and }sa\in bR\}.$$
\bpr
Let $R$ be a ring and $\mathbb L(R)$ be the category of principal left ideals. Let $\rho(a,s,b): Ra\rightarrow Rb$ be a morphism in $\mathbb L(R)$. Then 
\begin{enumerate}
	\item $\rho(a,s,b)$ is epimorphism if and only if $as \mc L b$
	\item $\rho(a,s,b)$ is a split monomorphism if and only if $a \mc R as$
	\item $\rho(a,s,b)$ is an isomorphism if and only if $a \mc R as \mc L b$
\end{enumerate}
\epr
\bp
$\rho(a,s,b): Ra\rightarrow Rb$ such that for all $x\in Ra$, $\rho(a,s,b): x\rightarrow xs$ where $as\in Rb$. So it is easy to observe that $\rho(a,s,b)$ is epimorphism if and only if $Ras = Rb$ and so $as \mc L b$. Now $\rho(a,s,b)$ is a split monomorphism if and only if there is a $\sigma = \rho(b,s',a): Rb\rightarrow Ra$ such that $\rho\sigma = I_{Ra}$ which implies $ass' = a$ and so $a \mc R as$. $(3)$ follows from $(1)$ and $(2)$.
\ep
\bl Let $R$ be a ring. Then category $\mathbb L(R)$ is a proper category.
\el
\bp
 If $\rho(a,s,b) = \rho(a,1,b)$ where $a\cdot 1 \in Rb$, then $\rho(a,s,b) = j(Ra,Rb)$ is an inclusion. Identity mapping and set inclusions of principal left ideals are morphisms in the category. So  $\mathbb L(R)$ is a category with subobjects. $\{\rho(a,1,b): a\in Rb\}$ is a choice of subobjects in the category $\mathbb L(R).$ If $\rho(a,s,b)$ is a morphism in $\mathbb L(R)$, then im $\rho(a,s,b)$ = $R{as}$ and $\rho(a,s,b) = \rho(a,s,as)\cdot \rho(as,1,b)$ gives the canonical  factorization of $\rho(a,s,b)$ in $\mathbb L(R)$. Let $\rho(a,1,b): Ra\subseteq Rb$ be the inclusion which splits. If inclusion splits, canonical factorization is unique(\cite{KSS}). Hence $\mathbb L(R)$ is a category with subobjects, every inclusion splits and every morphism has unique canonical factorization.\\
\i Now we define $\rho^d: \tt{v}\mathbb L(R) \rightarrow  Rd$
as $\rho^d(Ra) = \rho(a,s,d): Ra\rightarrow Rd$, where $Ra\in \tt v\mathbb L(R)$ and $as\in Rd$. Obviously $\rho^d(Ra) \in \mathbb L(R)(Ra,Rd)$ is well defined.\\If $Ra\subseteq Rb$, 
\begin{align*}
j(Ra,Rb)\rho^d(Rb) &= \rho(a,1,b)\cdot \rho(b,v,d) \text{ where } bv = qd \in Rd.\\
& = \rho(a,v,d)\\
& = \rho^d(Ra)
\end{align*}
Since $Ra\subseteq Rb$, $a = rb$ for some $r\in R$ and so $av = rbv = rqd \in Rd$. Hence $\rho(a,v,d) = \rho^d(Ra)$ and $j(Ra,Rb)\cdot \rho^d(Rb) = \rho^d(Ra)$. Hence $\rho^d$ is a cone in $\mathbb L(R)$. Since $\rho^d(Ra)$ has canonical factorization as $\rho^d(Ra) = \rho(a,s,d) = \rho(a,s,as)\cdot\rho(as,1,d)$, where $as\in Rd$. There exists some $s$ such that $Ras = Rd$. Thus $\rho^d$ is a proper cone in $\mathbb L(R)$ with vertex $Rd$ and hence $\mathbb L(R)$ is a proper category.
\ep 
Dually $\mathbb R(R)$ is also a proper category.
\bt 
The category $\mathbb L(R) [\mathbb R(R)]$ of principal left [right] ideals of a ring $R$ with right [left] translations as morphisms is a preadditive proper category.
\et
\bp
For $\rho(a,s,b), \rho(a,t,b)\in \mathbb L(R)(Ra,Rb)$, the addition is defined by 
$$\rho(a,s,b) + \rho(a,t,b) = \rho(a,s+t,b).$$ 
Under this addition homset\index{homset} $\mathbb L(R)(Ra,Rb)$ is an abelian group, further if restricted to $\mathbb L(R)(Ra,Ra)$, it is a ring.	\\
\i For $\rho(a,s,b), \rho(a,t,b), \rho(b,u,c),\rho(b,v,c) \in \mathbb L(R)$,
\begin{align*}
\rho(a,s,b)[\rho(b,u,c) + \rho(b,v,c)] &= \rho(a,s,b)[\rho(b,u + v,c)]\\
& = \rho(a,s(u + v),c) = \rho(a,su + sv,c)\\
& = \rho(a,su,c) + \rho(a,sv,c)\\
& = \rho(a,s,b)\rho(b,u,c) + \rho(a,s,b)\rho(b,v,c).
\end{align*}
Thus the composition of morphisms is bilinear. Zero ideal is the \tt {zero object} in $\mathbb L(R)$. Hence for a ring $R$, the categories $\mathbb L(R)$ and $\mathbb R(R)$ are preadditive proper categories. 
\ep 

\i Let $R$ be a ring with the property that any set of principal ideals which is bounded above has a maximal element and principal ideals of ring form a relatively complemented lattice with respect to the partial order induced from strict preorder. The join and meet are defined by; $Ra \vee Rb = Ra + Rb$ and $Ra \wedge Rb = Ra \cap Rb$ where $Ra + Rb$ is the smallest principal ideal containing both $Ra$ and $Rb$. The trivial ideals $R$ and $(0)$ are the bounds. Then $\mathbb L(R)$ and dually $\mathbb R(R)$ are $RR-$proper categories.  It is easy to see that the ideal categories of Euclidean domains are $RR-$proper categories. \\
\i In the $RR-$proper category $\mathbb L(R)$, by Theorem 3, the set of all proper cones $\{\rho^a:\tt v \mathbb L(R)\rightarrow Ra, \text{ for all }Ra\in \tt v\mathbb L(R)\}$ is a ring with the binary operations given below; for $Ra,Rb,Rc \in \tt v\mathbb L(R)$, 
$$(\rho^a\cdot\rho^b)(Rd) = \rho^a(Rd)\cdot(\rho^b(Ra))^o$$
$$\rho^a+\rho^b = (\rho^a\oplus\rho^b)^\ast$$
where $(\rho^a\oplus\rho^b)^\ast(Rd) = (\rho^a\oplus\rho^b)(Rd) e(Ra\vee Rb,d_0)$,
$$ (\rho^a\oplus\rho^b)(Rd) = \rho^a(Rd)j(Ra,Ra\vee Rb) + \rho^b(Rd)j(Rb,Ra\vee Rb)$$ and $d_0 = max\{im (\rho^a\oplus\rho^b)(Rd), |Rd\in \tt v\mathbb L(R)\}$. We denote this ring of proper cones by $\mc P\mathbb L(R)$. 

\i By Theorem 4, the left(right) ideal category of the ring $\mc P\mathbb L(R)$ are preadditive proper categories $\mathbb L(\mc P\mathbb L(R))$($\mathbb R(\mc P\mathbb L(R))$). Also the left(right) ideals form a relatively complemented lattice with meet intersection and join as the smallest principal ideal cotaining both. So $\mathbb L(\mc P\mathbb L(R))$ and $\mathbb R(\mc P\mathbb L(R))$ are $RR-$proper categories and hence the set of all proper cones $\mc P(\mathbb L(\mc P\mathbb L(R)))$ and $\mc P(\mathbb R(\mc P\mathbb L(R)))$ are rings.

\end{document}